\documentclass[a4paper,abstracton]{scrartcl}
\usepackage[utf8]{inputenc}
\usepackage{amsfonts,amsthm,amssymb,amsmath}
\usepackage[inline]{enumitem}
\usepackage{color}
\usepackage{graphicx} 
\usepackage{authblk}

\newtheorem{thm}{Theorem}
\newtheorem{lem}[thm]{Lemma}

\theoremstyle{remark}

\newtheorem{clm}{Claim}

\definecolor{pink}{RGB}{219, 48, 122}

\DeclareMathOperator{\AC}{\bf{AC}}

\title{The role of the Axiom of Choice in proper and distinguishing colourings}           
{}                 

\author{Marcin Stawiski\\stawiski@agh.edu.pl}            
\affil{AGH University,\\ Faculty of Applied Mathematics, \protect\\al. Mickiewicza 30, 30-059 Krakow, Poland}  
{}                     
{}

\begin{document}
\maketitle


\begin{abstract}
Call a colouring of a graph \emph{distinguishing} if the only automorphism  which preserves it is the identity. We investigate the role of the Axiom of Choice in the existence of certain proper or distinguishing colourings in both vertex and edge variants with special emphasis on locally finite connected graphs. We show that  every locally finite connected graph has a distinguishing colouring with at most countable number of colours or every locally finite connected graph has a proper colouring with at most countable number of colours if and only if Kőnig's Lemma holds. This statement holds for both vertex and edge colourings. Furthermore, we show that it is not provable in ZF that such colourings exist even for every connected graph with maximum degree 3.
We also formulate a few conditions about distinguishing and proper colourings which are equivalent to the Axiom of Choice.
   \end{abstract}


\section{Introduction}
 Let $c$ be a vertex or an edge colouring of a graph $G$. We say that an automorphism $\varphi$ of  $G$ \emph{preserves} $c$ if each vertex of $G$ is mapped to a vertex of the same colour or each edge of $G$ is mapped to an edge of the same colour.
Call a colouring $c$  \emph{distinguishing} if the only automorphism which preserves $c$ is the identity. If a colouring $c$ is a mapping into ordinal numbers (or any well-ordered set) we can think about the number of colours in $c$.  The \emph{distinguishing number} $D(G)$ of a graph $G$ is the least number of colours in a distinguishing vertex colouring of $G$. Similarly, the \emph{distinguishing index} $D'(G)$ of a graph $G$ is the least number of colours in a distinguishing edge colouring of $G$.
Distinguishing vertex colourings were introduced by Babai \cite{BAB} in 1977 under the name \emph{asymmetric} colourings during his study of the complexity of the graph isomorphism problem \cite{babaiisomorphism}. Distinguishing edge colourings were introduced by Kalinowski and Pil\'sniak \cite{KP}.

In this paper we study proper and distinguishing colourings in ZF, hence without assuming the Axiom of Choice. Proper vertex colourings in ZF were investigated by Galvin and Komjáth \cite{galvin}. They proved that the existence of the chromatic number of each graph is equivalent to the Axiom of Choice. Distinguishing colourings and proper edge colourings in ZF were not previously investigated. 

In most of the papers about infinite graphs some version of the Axiom of Choice is used though  not always explicitly. The most popular methods often involve Zorn's Lemma or Kőnig's Lemma. In particular, proofs of general bounds by a function of $\Delta(G)$ for chromatic number \cite{bruijn}, chromatic index \cite{lake}, distinguishing number \cite{LPS} and distinguishing index \cite{PS} of connected infinite graphs  all use Kőnig's Lemma in the case of locally finite graphs and the Axiom of Choice in the form of Hessenberg's Theorem for general bounds. We show that in all these cases the use of  Kőnig's Lemma or respectively the Axiom of Choice is necessary.

Similar problems for graphs without the assumption of connectivity were previously investigated for proper vertex colourings. The statement that every graph has a proper vertex colouring using at most two colours  if and only if each of its finite subgraphs has such a  colouring is equivalent to the Axiom of Choice for Pairs. If we replace two colours with three colours, then we obtain the statement equivalent to the Prime Ideal Theorem. See \cite{Herrlich} p. 109--116 for details and further examples.

Arguably, most of the results related to  proper or distinguishing colourings in graph theory concern only locally finite connected graphs.
From the results mentioned in the previous paragraph, it follows that one cannot prove in ZF that every locally finite connected  graph has a distinguishing or a proper colouring with at most countable number of colours. We show that one cannot prove the existence of such colourings in ZF even in the simplest case of connected graphs with maximum degree 3.

\section{Preliminaries}

By a \emph{cardinal number} we mean an \emph{initial ordinal} i.e. an ordinal which is not equinumerous with any smaller ordinal. For every set there exists a cardinal number equinumerous with it if and only if the Axiom of Choice holds. 

 Well-Ordering Theorem states that for every set $X$ there exists a well-order  on $X$. Well-Ordering Theorem is  equivalent to the Axiom of Choice.

We now present some weak choice principles. The axiom $\AC^\omega_{\leqslant \kappa}$ states that every countable family of non-empty sets of cardinality at most $\kappa$ has a choice function. The axiom $\AC^\omega_{fin}$ states that every countable family of non-empty finite sets has a choice function. The axiom $\AC^\omega_{2}$ is the same as $\AC^\omega_{\leqslant 2}$. The axiom $\AC^\omega_{fin}$ is equivalent to Kőnig's Lemma stating that every locally finite infinite connected graph has a ray. The axiom $\AC^\omega_{fin}$ is also equivalent to the statement that every countable union of finite sets is countable. More about the Axiom of Choice, weak choice principles and their equivalent forms may be found in the extensive monograph of Howard and Rubin \cite{consequences}.

Let $G$ be a graph. Denote by $\Delta(G)$ the supremum over the
 degrees of all vertices of $G$. If there exists a vertex $v \in V(G)$ such that $d(v)=\Delta(G)$, then $\Delta(G)$ is called the \emph{maximum degree} of $G$. Graphs with maximum degree 3 are called \emph{subcubic}. We say that a graph is \emph{locally finite} if each of its vertices has finite degree.

Let $\Gamma$ be a group acting on a set $\Omega$ and let $A$ be a subset of $\Omega$. The \emph{orbit} of  $A$ is the set $\{\varphi(a): a\in A, \varphi \in \Gamma \}$. We say that $A$ is \emph{fixed} if every $\varphi \in \Gamma$ acts trivially on $A$ i.e. if $\varphi(a)=a$ for every $a\in A$.  We say that $A$ is \emph{stabilized} if for every $\varphi \in \Gamma$, we have $\varphi(A) \subseteq A$.  In the definitions in this paragraph, if $A=\{ a \}$ is a singleton, then we often refer to $a$ instead of $\{ a\}$. 
An \emph{automorphism} of a graph $G$ is a bijection $\varphi: V(G) \rightarrow V(G)$ such that $uv$ is an edge in $G$ if and only if $\varphi(u)\varphi(v)$ is an edge in $G$. They form a group with composition as the operation. If not written explicitly, the meaning of $\Gamma$ and $\Omega$ shall follow from the context. In this paper $\Gamma$ is usually a group of some automorphisms of a graph $G$ and $\Omega$ is a set of some vertices of $G$ or some edges of $G$. 

Colourings in this paper are not necessarily proper unless stated otherwise.
For notions which are not defined here, see \cite{diestel} or \cite{Jech}.

\section{The Axiom of Choice in proper and distinguishing colourings}\label{chapter:ac}

Let $\kappa$ be an arbitrary non-zero cardinal. Call a family $\mathcal{A}=\{A_i : i \in \omega \}$ \emph{acceptable} if $\mathcal{A}$ is a countable family of pairwise disjoint non-empty sets. We say that a family $\mathcal{A}$ is \emph{almost $\kappa$-acceptable} if $\mathcal{A}$ is acceptable and every set in $\mathcal{A}$ has cardinality less than $\kappa$. We say that  $\mathcal{A}$ is \emph{$\kappa$-acceptable} if it is almost $\kappa^+$-acceptable. In other words, $\mathcal{A}$ is \emph{$\kappa$-acceptable} if  $\mathcal{A}$ is acceptable and every set in $\mathcal{A}$ has cardinality at most $\kappa$.

Let $\mathcal{A}=\{A_i : i \in \omega \}$ be an acceptable family and let $Y=\bigcup \mathcal{A}$. Let $Z=\{ z_i:i \in \omega\}$ and $Z'=\{ z_i':i \in \omega\}$ be disjoint sets which are also disjoint from $Y$. We now define graphs $G_\mathcal{A}$ and $H_\mathcal{A}$ by
\begin{align*}
V(G_\mathcal{A})&=V(H_\mathcal{A})=Y \cup Z \cup Z', \\
E(G_\mathcal{A})&=\{ z_i z_{i+1}:i\in \omega\} \cup \{z_i z_i' :i\in \omega\} \cup \{a_i z_i' : i \in \omega, a_i \in A_i \}  ,\\
E(H_\mathcal{A})&=E(G_\mathcal{A}) \cup \{ab:a\neq b, a,b \in A_i, i \in \omega \}.
\end{align*}
 From the definitions of $G_\mathcal{A}$ and $H_\mathcal{A}$ it follows  that every vertex in $Z \setminus \{z_0\}$ has degree 3 in both graphs, and vertex $z_0$ has degree 2. For the rest of the paragraph, assume that for every $i\in \omega$ the set $A_i$ is well-orderable. The vertex $z'_i$ has degree $|A_i|+1$ for every $i\in \omega$ in both $G_\mathcal{A}$ and $H_\mathcal{A}$. Every vertex in $Y$ has degree 1  in $G_\mathcal{A}$. However, every vertex $a \in A_i$ has degree $|A_i|$ in $H_\mathcal{A}$ since it has edges to every other vertex in $A_i$ and to $z'_i$. Summarizing, we obtain $\Delta(G_\mathcal{A})=\Delta(H_\mathcal{A})=\max \{ 3 , \sup\{ | A_i|+1: i \in \omega \}\}$.

\begin{clm}\label{obserwacja}
Let $\mathcal{A}=\{A_i :i \in \omega \}$ be an acceptable family. Then for every natural number $i$, the vertices in $A_i$ form an orbit with respect to the groups of automorphisms of $G_\mathcal{A}$ and $H_\mathcal{A}$. The rest of the vertices in both graphs are fixed with respect to these groups.
\end{clm}
\begin{proof}

First we show that $z_0$ is fixed in both graphs. Suppose that $z_0$ may be mapped into $z_i$ for some $i \neq 0$. Let $R_i$ be the maximal induced ray with endvertex $z_i$. Clearly $R_0$ is a tail of the ray $R$ induced by $Z$. Let $P$ be a maximal induced path with endvertex $z_0$ which is edge disjoint from $R$, and let $P_i$ be a maximal induced path with endvertex $z_i$ which is edge disjoint from $R_i$ and contains $z_0$.  If $P_i$ contains $z_0$, then the length of $P_i$ is larger than the length of $P$. Notice that in this case $P_i$ is the longest induced path with endvertex $z_i$ which is edge disjoint from $R_i$. This leads to contradiction because $P_i$ cannot be mapped into $P$. Hence, $z_0$ is fixed. The ray $R$ is the only induced ray with endvertex $z_0$. Therefore, $R$ is fixed.
\begin{figure}\centering\includegraphics[scale=1]{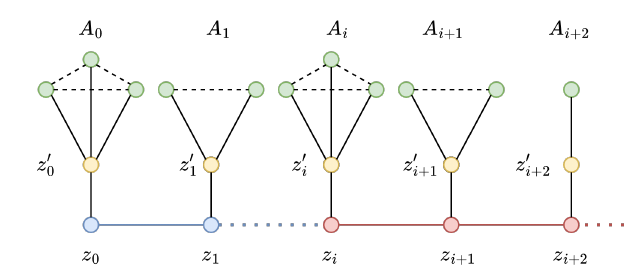}\caption{Graphs $G_{\mathcal{A}}$ and $H_\mathcal{A}$. Graph $H_\mathcal{A}$ is obtained by adding the dashed edges. The ray $R_i$ is pink, and the remaining part of the ray $R$ is blue.}\end{figure}

Since every vertex of the fixed set $Z$ has exactly one neighbour outside $Z$, all of these neighbours are fixed. Hence, $Z'$ is fixed, and $A_i$ is stabilized for every $i \in \omega$. The vertices in $A_i$ form an independent set in $G_\mathcal{A}$ (or a clique in $H_\mathcal{A}$). Hence, they form an orbit with respect to the group of automorphisms of $\mathcal{G}_A$, and also with respect to the group of automorphisms of $\mathcal{H}_A$.
\end{proof}

Notice that from Claim \ref{obserwacja} it follows that the group of automorphisms of $G_\mathcal{A}$ is the same as the group of automorphism of $H_\mathcal{A}$. Now, we prove a lemma which allow us to restrict part of the later considerations to the problem of the existence of the distinguishing number for graphs of the form $G_\mathcal{A}$. With the lemma below we are able to simultaneously obtain results about distinguishing colourings and proper colourings in both vertex and edge versions.

\begin{lem}\label{lem:rownowazneliczbie} Let $\mathcal{A}$ be an acceptable family. Then the following conditions are equivalent.
\begin{enumerate}[labelindent=\parindent,leftmargin=*, align=left,label = \textnormal{\alph*)}]
    \item\label{NWSR:11} There exists the distinguishing number of $G_\mathcal{A}$.
    \item\label{NWSR:12} There exists the distinguishing index of $G_\mathcal{A}$.
    \item\label{NWSR:13} There exists the chromatic index of $G_\mathcal{A}$.
    \item\label{NWSR:14} There exists the chromatic number of $H_\mathcal{A}$.
\end{enumerate}
\end{lem}
\begin{proof}
By Claim \ref{obserwacja} if $c$ is a distinguishing vertex colouring of $G_\mathcal{A}$, then for every $i \in \omega$ vertices in $A_i$ have distinct colours. If for every $i\in \omega$ and each vertex $v\in A_i$ we colour the edge $vz_i'$ with colour $c(v)$, and we colour the rest of the edges of $G_\mathcal{A}$ arbitrarily, then we obtain a distinguishing edge colouring of $G_\mathcal{A}$. Hence, condition \ref{NWSR:11} implies \ref{NWSR:12}.

Now, let $c$ be a distinguishing edge colouring of $G_\mathcal{A}$. Let $c'$ be a colouring in which  the edges incident to vertices in $Y$ have the same colour as in $c$, the edges between $Z$ and $Z'$ are coloured with the same new colour, and the edges between the vertices in $Z$ are coloured alternately with two new colours. The colouring $c'$ is a proper edge colouring. Therefore, condition \ref{NWSR:12} implies condition \ref{NWSR:13}.

Let $c$ be a proper edge colouring of $G_\mathcal{A}$. We  now define a proper vertex colouring $c'$ of $H_\mathcal{A}$. First, for every $i\in \omega$, we colour each vertex $v\in A_i$ with colour $c'(v)=c(vz_i')$. Next, we colour the vertices in $Z'$ with the same new colour, and we colour the vertices in $Z$ alternately with two new colours. Colouring $c'$ is a proper vertex colouring of $H_\mathcal{A}$. Hence, condition \ref{NWSR:13} implies condition \ref{NWSR:14}.

Implication between \ref{NWSR:14} and \ref{NWSR:11} follows directly from Claim \ref{obserwacja} since every proper vertex colouring of $H_\mathcal{A}$ is a distinguishing vertex colouring of $G_\mathcal{A}$.
\end{proof}

We can now proceed to the study of relations between the existence of certain colourings and the Axiom of Choice. The first step is Lemma \ref{lem:istnieniefunkcjiwyboru}, which shows that the existence of the distinguishing number of $G_\mathcal{A}$ implies the existence of a choice function for $\mathcal{A}$.

\begin{lem}\label{lem:istnieniefunkcjiwyboru} Let $\mathcal{A}$ be an acceptable family and assume that there exists the distinguishing number of $G_\mathcal{A}$. Then there exists a choice function for $\mathcal{A}$.
\end{lem}
\begin{proof}
Let $c$ be a vertex colouring of $G_\mathcal{A}$ with elements of some cardinal number $\kappa$.  Then $f(A)= \arg \min \{ c(a) : a\in A\}$ is a choice function for $\mathcal{A}$.
\end{proof}

We now prove the next lemma which in the case of non-zero natural number $k$ and $k$-acceptable family $\mathcal{A}$ allows us to construct a distinguishing colouring of $G_\mathcal{A}$ using a choice function for $\mathcal{A}$.

\begin{lem} \label{choice_0} Let $k$ be an arbitrary non-zero natural number and assume $\AC^\omega_{\leqslant k}$. Then for every $k$-acceptable family $\mathcal{A}$ graph $G_\mathcal{A}$ has distinguishing number at most $k$.
\end{lem}
\begin{proof}
 The proof is by induction on $k$. Let $\mathcal{A}$ be a $k$-acceptable family. If $k=1$, then $G$ has no non-trivial automorphism. Hence, its distinguishing number is equal to 1. Assume that $k\geqslant 2$, and that the statement of the lemma holds for every $l < k$. Let $f$ be a choice function for $\mathcal{A}$. From the inductive hypothesis $G_\mathcal{A}-f(\mathcal{A})$ has a distinguishing vertex colouring $c'$ using at most $k-1$ colours. Colouring $c$ which agrees with colouring $c'$ on $V(G_\mathcal{A})\setminus f(\mathcal{A})$ and which assigns the rest of vertices of $G_\mathcal{A}$ the same new colour is a distinguishing colouring using at most $k$ colours. \end{proof}

Lemmas \ref{lem:rownowazneliczbie}--\ref{choice_0} allows us to formulate the following corollary  about the existence of certain parameters for $k$-acceptable families in the case of finite $k$.

\begin{thm}\label{tw:knaturalne}
Let $k\geqslant 2$ be an arbitrary natural number. Then the following conditions are equivalent.
\begin{enumerate}[labelindent=\parindent,leftmargin=*, align=left,label = \textnormal{\alph*)}]
    \item $\AC^\omega_{\leqslant k}$. \label{itm:knaturalne1}
    \item For every $k$-acceptable family $\mathcal{A}$ the graph $G_\mathcal{A}$ has the distinguishing number.
    \item For every $k$-acceptable family $\mathcal{A}$ the graph $G_\mathcal{A}$ has the distinguishing index. 
    \item For every $k$-acceptable family $\mathcal{A}$ the graph $G_\mathcal{A}$ has the chromatic index. 
    \item For every $k$-acceptable family $\mathcal{A}$ the graph $H_\mathcal{A}$ has the chromatic number.
\end{enumerate}
\end{thm}

In particular for $k=2$ condition \ref{itm:knaturalne1} is the axiom  $\AC^{\omega}_2$ which is independent of ZF. It follows that in ZF one cannot prove the existence of the above parameters even for every connected subcubic graph.

Theorem \ref{tw:knaturalne} tells us about the existence of certain parameters for connected graphs with finite maximal degree. Now, we  establish the relations between  Kőnig's Lemma and the existence of proper colourings and distinguishing colourings  using at most countable number of colours in the case of locally finite connected graphs.

\begin{thm}
The following conditions are equivalent.
\begin{enumerate}[labelindent=\parindent,leftmargin=*, align=left, label = \textnormal{(KL\arabic*)}, resume=koniglemma]
    \item [\textnormal{(KL)}] Kőnig's Lemma.\label{nwsr:1}
    \item Every infinite  locally finite connected graph has the distinguishing number.\label{nwsr:2}
    \item Every infinite  locally finite connected graph has the distinguishing index.\label{nwsr:3}
    \item Every infinite  locally finite connected graph has the chromatic index.\label{nwsr:indeks1}
    \item Every infinite  locally finite connected graph has the chromatic number.\label{nwsr:liczba1}
    \item For every almost $\aleph_0$-acceptable family $\mathcal{A}$ the graph $G_{\mathcal{A}}$ has the  distinguishing number.\label{nwsr:4}
    \item For every almost $\aleph_0$-acceptable family $\mathcal{A}$ the graph $G_{\mathcal{A}}$ has the  distinguishing index.\label{nwsr:5} 
    \item For every almost $\aleph_0$-acceptable family $\mathcal{A}$ the graph $G_{\mathcal{A}}$ has the  distinguishing index.\label{nwsr:indeks2}  
    \item For every almost $\aleph_0$-acceptable family $\mathcal{A}$ the graph $H_{\mathcal{A}}$ has the  chromatic number.\label{nwsr:liczba2}
\end{enumerate}
\end{thm}
\begin{proof}
First, we show that Kőnig's Lemma implies conditions \ref{nwsr:2}--\ref{nwsr:liczba1}. Let $G$ be an infinite locally finite connected graph.  Let $v$ be some vertex of $G$. For a natural number $d$ denote $B(v,d)=\{x\in V(G): d(v,x)=d \}$. By local finiteness of $G$ each $B(v,d)$ is finite.  By connectivity of $G$ the vertex set of $G$ may be represented as the countable union of finite sets $V(G)=\bigcup\{B(v,d):d < \omega \}$. Recall that Kőnig's Lemma is equivalent to the statement that the sum of every countable family of finite sets is countable. Hence, $V(G)$ is countable.
As $E(G) \subseteq V(G) \times V(G)$, then $E(G)$ is also  countable. Since both  sets $V(G)$ and $E(G)$ are countable, we can obtain the desired colourings by assigning to each vertex (edge respectively) a unique natural number.

Implications $\ref{nwsr:2}\Rightarrow \ref{nwsr:4}$, $\ref{nwsr:3}\Rightarrow \ref{nwsr:5}$, $\ref{nwsr:indeks1}\Rightarrow \ref{nwsr:indeks2}$ and $\ref{nwsr:liczba1}\Rightarrow \ref{nwsr:liczba2}$ are trivial. The equivalence of the conditions \ref{nwsr:4}--\ref{nwsr:liczba2} follows from Lemma \ref{lem:rownowazneliczbie}. 

It remains to show that the condition \ref{nwsr:4} implies Kőnig's Lemma. From \ref{nwsr:4} and Lemma
\ref{lem:istnieniefunkcjiwyboru}, we have that for every countable family of finite sets there exists its choice function. This is the axiom $\AC^\omega_{fin}$, which is  equivalent to Kőnig's Lemma. \end{proof}

As we have shown, the existence of the distinguishing number of $G_\mathcal{A}$ for every almost $\aleph_0$-acceptable family $\mathcal{A}$ is equivalent to Kőnig's Lemma and therefore to the Axiom of Countable Choice for Finite Sets. One may think that the existence of the distinguishing number of every graph of the form $G_\mathcal{A}$ for some acceptable family $\mathcal{A}$ is equivalent to the Axiom of Countable Choice. It turns out that this condition is much stronger and it implies the full Axiom of Choice.

\begin{thm}\label{liczbaimplikujeac}
If for every acceptable family $\mathcal{A}$ the graph $G_{\mathcal{A}}$ has the distinguishing number, then the Axiom of Choice holds.
\end{thm}
\begin{proof}
Let $X$ be a non-empty set and let $\mathcal{A}$ be an acceptable family such that $X \in \mathcal{A}$. By the assumption  there exists a distinguishing vertex colouring $c$ of the graph $G_\mathcal{A}$ using colours from some cardinal $\kappa$. As the colouring $c$ is distinguishing, the elements of $X$ have distinct colours in $c$. It follows that $c|_X$ is an injection from $X$ to cardinal number $\kappa$. Hence,  the Well-Ordering Theorem holds and so does the Axiom of Choice.
\end{proof}

Theorem \ref{liczbaimplikujeac} allows to formulate a list of conditions equivalent to the Axiom of Choice. The conditions \ref{nwsrAC:1j}--\ref{nwsrAC:4j} in the theorem below are equivalent to their restrictions to connected graphs. Recall that the equivalence of the Axiom of Choice, the existence of the chromatic number of every graph, and the existence of the chromatic number of every connected graph was proved by Galvin and Komjáth \cite{galvin}. 

\begin{thm}\label{tw:kolorowaniarownowazneac}
The following conditions are equivalent.
\begin{enumerate}[labelindent=\parindent,leftmargin=*, align=left, label = \textnormal{(AC\arabic*)}, series=ac]
    \item [\textnormal{(AC)}] The Axiom of Choice.\label{nwsrAC:0}
    \item Every graph has  the distinguishing number. \label{nwsrAC:1j}
    \item Every  graph without a component isomorphic to $K_1$ or $K_2$ has  the distinguishing index.\label{nwsrAC:2j}
    \item Every  graph has the chromatic index.\label{nwsrAC:3j}
    \item Every  graph has the chromatic number. \label{nwsrAC:4j}
    \item For every acceptable family $\mathcal{A}$ the graph $G_{\mathcal{A}}$ has  the distinguishing number.\label{nwsrAC:1i}
    \item For every acceptable family $\mathcal{A}$ the graph $G_{\mathcal{A}}$ has   the distinguishing index.\label{nwsrAC:2i}
    \item For every acceptable family $\mathcal{A}$ the graph $G_{\mathcal{A}}$ has   the chromatic index.\label{nwsrAC:3i}
    \item For every acceptable family $\mathcal{A}$ the graph $H_{\mathcal{A}}$ has   the chromatic number.\label{nwsrAC:4i}
\end{enumerate}
\end{thm}
\begin{proof}
From the Well-Ordering Theorem we can well-order the set of vertices and the set of edges of a given graph and then colour each vertex (edge respectively) of the said graph with a unique colour. This means that the Axiom of Choice implies conditions \ref{nwsrAC:1j}--\ref{nwsrAC:4j}. Each of the condition \ref{nwsrAC:1j}--\ref{nwsrAC:4j} implies its restriction to connected graphs and also the corresponding condition \ref{nwsrAC:1i}--\ref{nwsrAC:4i}.
By Lemma \ref{lem:rownowazneliczbie} conditions \ref{nwsrAC:1i}--\ref{nwsrAC:4i} are equivalent. The implication between condition \ref{nwsrAC:1i} and the Axiom of Choice is Theorem \ref{liczbaimplikujeac}.
\end{proof}

\end{document}